\documentclass{sn-jnl}

\usepackage{graphicx}%
\usepackage{multirow}%
\usepackage{amsmath,amssymb,amsfonts}%
\usepackage{amsthm}%
\usepackage{mathrsfs}%
\usepackage[title]{appendix}%
\usepackage{xcolor}%
\usepackage{textcomp}%
\usepackage{manyfoot}%
\usepackage{booktabs}%
\usepackage{algorithm}%
\usepackage{algorithmicx}%
\usepackage{algpseudocode}%
\usepackage{listings}%
\usepackage{comment}

\def\eqindistn{\mathop{=^{\mkern-16mu D}\,}}

\begin{document}

\title{Quantile Fourier regressions for decision making under uncertainty}

\author[1]{\fnm{Arash} \sur{Khojaste}}\email{akhojaste@umass.edu}

\author*[2]{\fnm{Geoffrey} \sur{Pritchard}}\email{g.pritchard@auckland.ac.nz}

\author[1]{\fnm{Golbon} \sur{Zakeri}}\email{gzakeri@umass.edu}

\affil[1]{\orgdiv{Dept. of Mechanical and Industrial Engineering}, \orgname{University of Massachusetts -- Amherst}, \orgaddress{\city{Amherst}, \postcode{01002}, \state{MA}, \country{United States}}}

\affil*[2]{\orgdiv{Dept. of Statistics}, \orgname{University of Auckland}, \orgaddress{\street{Private Bag 92019}, \city{Auckland}, \postcode{1142}, \country{New Zealand}}}


\abstract{We consider Markov decision processes arising from a Markov model of an underlying natural
phenomenon. Such phenomena are usually periodic ({\it e.g.} annual) in time, and so the Markov processes
modelling them must be time-inhomogeneous, with cyclostationary rather than stationary behaviour.
We describe a technique for constructing such processes that allows for periodic variations both
in the values taken by the process and in the serial dependence structure. We include two illustrative
numerical examples: a hydropower scheduling problem and a model of offshore wind power integration.}


\maketitle
\section{Introduction}

We consider stochastic optimization problems in which the underlying randomness takes the form of a
discrete-time univariate stochastic process with periodically varying behaviour.
Our approach will be to approximate the underlying stochastic process by a finite-state Markov chain
in which both the interpretation of the states and the Markov transition probabilities are
periodically varying (time-inhomogeneous).
The optimization problem is then modelled as a Markov decision process.

Many real-world planning problems involve decisions of two kinds, which we characterize as ``investment'' and ``operations''.
In the investment phase, decisions are taken which create or preclude possible courses of action in the later operations phase.
Examples include capacity planning for electric power, management and operations in forestry, and literal financial
investments. The early decisions ({\it e.g.} on the construction or retirement of power-generation plants) can have a large effect
on the cost and robustness of later operations ({\it e.g.} efficient and reliable dispatch of electricity).
Uncertainty plays a key role: the investment decisions must be made without full knowledge of the operational environment
in which their consequences will play out.


The operations phase typically represents a long-term future, and so it is often further subdivided into many stages.
At each stage, some of the uncertainty is resolved, and recourse actions become possible in response to the new information.
This is the paradigm of multi-stage stochastic programming (see {\it e.g.} \cite{AhmedKingParija,DingHuBie,DingQuShahidehpour,SabetYazdani}).
The uncertainty in such problems can be represented by scenario trees that branch at each operational stage.
For example, investment decisions for electricity generation would be followed by a multitude of hour-long operational stages
(8760 of them in each year, perhaps for many future years), with multiple random quantities realized at each stage.
It is evident that due to the exponential growth in the number of scenarios and the distant time horizon,
such problems can become prohibitively large. Much research effort has been spent on developing decomposition approaches
and approximation schemes that iteratively solve a limited version of the problem in pursuit of a solution
(see {\it e.g.} \cite{AhmedKingParija,KleywegtShapiroH,ZakeriPhilpottRyan}). 


This paper takes a different approach: we propose to represent the operations phase via the stationary behaviour
of a discrete Markov decision process (MDP). While the usual exposition of MDPs posits a single endlessly repeated problem
(given that the system is in one of finitely many possible states, choose one of finitely many possible actions which will
determine the probability distribution of the state in the next incarnation), it is readily adaptable to models
in which the states, actions, and associated parameters vary in a periodic cycle. The solution is then a cyclostationary
rather than a stationary process.

This succinct representation of uncertainty can then be embedded within the investment phase of the problem, allowing
the effects of investment decision-making on subsequent operations to be easily evaluated.
The present paper, however, will confine its focus to the modelling and solution of the operations phase.
We will outline a systematic approach for constructing a suitable Markov decision process via quantile Fourier regression.
Quantile regressions are ideal for modelling continuously and periodically varying random phenomena (\cite{Pritchard_hydro})
in which the probability distributions must be different at each point in the periodic cycle ({\it e.g.} on each day of the year).
We also describe an apparently new method of fitting time-inhomogeneous Markov transition probabilities,
so that not only the random values themselves, but also their serial dependence structure, can be subject to periodic variation.

We illustrate the techniques developed in this paper with two applications. In section \ref{s:hydro} we discuss a
hydro-thermal power scheduling problem with seasonal variation in inflows. In section \ref{s:offshore_wind} we consider a
more ambitious model of offshore windpower integration subject to both diurnal and annual variations.

\section{Quantile Fourier regression}
\label{s:QFR}

A discrete-time stochastic process $(X_t)_{t=0}^\infty$ is {\em cyclostationary} if, for any integers $m_1,\ldots, m_k$,
the joint probability distribution of $X_{t+m_1},\ldots,X_{t+m_k}$ is a periodic function of $t$.
For the purposes of the present paper, the period $T$ will always be an integer; thus
$$ (X_{t+m_1},\ldots,X_{t+m_k}) \eqindistn (X_{t+m_1+T},\ldots,X_{t+m_k+T}) .
  $$
We also assume that $T$ is known.
For further and more general discussion of cyclostationary processes, the reader is referred to \cite{GNP06}.

In modelling some phenomenon as a cyclostationary process, we might perhaps begin by describing the
probability distribution of $X_t$ as a function of $t$.
We shall try to avoid making specific assumptions about the shape(s) of this distribution.
However, if the period $T$ is a large integer---{\it e.g.} a process with hourly or daily
time steps and annual periodicity---it is desirable to assume some structure in the way the
distribution varies as a periodic function of $t$, rather than attempting to fit $T$ completely
unrelated distributions.
(The most serious objection to the latter approach is that each such distribution must be estimated
from a subsample comprising only \textcolor{black}{a fraction $\frac1T$}
of the available data; this subsample will usually
be too small to do the job well.)

\textcolor{black}{The present paper will achieve this end via a technique
that effectively treats time as a continuous variable:}
quantile Fourier regression.
For a fixed $p\in(0,1)$, the $p$-quantile of the distribution of $X_t$ is modelled as a
function $q_{p}(t)$ with period $T$ drawn from a finite-dimensional space of such functions:
\begin{equation}
\label{eq:QR}
P(X_t \leq q_{p}(t)) = p
  \qquad\hbox{ where }\qquad
  q_{p}(t) = \sum_{j=1}^d \beta_j b_j(t) .
\end{equation} 
Here $\beta_1,\ldots,\beta_d$ are coefficients to be fitted from data,
while $b_1(t),\ldots,b_d(t)$ are a basis for the chosen function space.
\textcolor{black}{(The coefficients $\beta_1,\ldots,\beta_d$,
and possibly also the basis functions $b_1(t),\ldots,b_d(t)$,
will be different for each $p$; we suppress this dependence in the
interest of notational brevity.)}
Suitable function spaces include periodic splines and wavelets
(\cite{Averbuch}). But in the present paper, we make the same choice
of basis as Jean-Baptiste Joseph Fourier: $d=2r+1$ is odd and
\begin{equation}
\label{eq:Fourier_basis}
  b_0(t) = 1 
  \quad\hbox{ and }\quad
  b_{2k-1}(t) = \cos(k\omega t),\quad  b_{2k}(t) = \sin(k\omega t),
  \hbox{ for } k=1,\ldots,r ,
\end{equation} 
where $\omega=2\pi/T$.
For small values of $r$, the quantile is thus constrained to a smooth and
fairly simple variation with $t$.

\textcolor{black}{
Estimation of the coefficients $\beta_j$ can be done by quantile regression,
a well-understood and computationally straightforward technique
(\cite{Koenker}).
Suppose that for each $k=1,\ldots,n$ we have observed a value $x_k$
at time $t_k$. The estimation is done by solving an optimization problem:
the coefficients $\beta_j$ should be chosen to minimize the value of
$$ \sum_{k=1}^n f_p\left(x_k - \sum_{j=1}^d \beta_j b_j(t_k)\right) ,
  $$
where $f_p$ denotes the function
$f_p(x) = \hbox{max}(p x, (p-1)x)$.
This optimization reduces to a linear programming problem.
For further details, the reader is referred to (\cite{Koenker}).}

After constructing such models for several different values of $p$,
we have a (partial) description of the distribution of $X_t$ throughout the
periodic cycle.

\textcolor{black}{A simple} special case of (\ref{eq:QR}) occurs when the
$b_j(t)$ are piecewise constant
\textcolor{black}{ ({\it i.e.} step functions) in $t$}.
This occurs, for example,
when there is annual periodicity and the stationary distribution
of $X_t$ is assumed to be constant within each month of the year, so that there are 12
(unrelated) distributions to fit. With this approach, estimation could hardly be simpler:
to estimate, say, the 75th percentile of the distribution for January, we need only
calculate the 75th percentile of all the data points observed in January of any year.
However, this is a poor choice of functional form when our ultimate intention is to incorporate the
resulting model into a stochastic optimization problem. The discontinuous changes in the
distribution at certain points in the periodic cycle are likely to distort optimal decisions
near those points (a ``month-end effect''). It is for this reason that the present paper
considers only continuous models for the quantiles.

It is possible that there is more than one periodic cycle. For example, a process with hourly
time steps may have both diurnal and annual periodicity. For simplicity, we assume in this paper that
the periods are commensurate, so that there is a single overall period which can be expressed
as an integer multiple of each individual period.
The modelling problem then reduces to finding a function space
for the quantiles which incorporates all of the periodic variations present.

\section{A Markov model of serial dependence}
\label{s:Markov_model}

Describing the probability distribution of $X_t$ for each $t$ does not
complete the modelling task; we should also consider the serial dependence
structure of the process.

\textcolor{black}{
The quantile representation of Section \ref{s:QFR} naturally reduces the univariate process $(X_t)$ to a finite-state one $(Z_t)$.
Given functions $q_{p_1}(t),\ldots,q_{p_{m-1}}(t)$ specifying quantiles
$p_1<\cdots<p_{m-1}$ of $X_t$, then we can define the random variable
$Z_t$ by $Z_t=i\iff q_{p_{i-1}}(t)\leq X_t < q_{p_i}(t)$
(``the system is in state $i$'') for $i=2,\ldots,m-1$.}
We also have extreme intervals $(-\infty, q_{p_1}(t))$
and $[q_{p_{m-1}}(t),\infty)$ represented by states $1$ and $m$
respectively. For convenience in referring to these, we define $p_0=0$,
$p_m=1$, $q_0(t)=-\infty$, and $q_1(t)=\infty$.

\textcolor{black}{
Perhaps the simplest model of serial dependence in a finite-state process
is the Markov chain. We are now well-placed to adopt such a model for $(Z_t)$.}

Inference of the Markov transition probabilities can be carried out in several
ways. The most straightforward is to fit a single transition
\textcolor{black}{probability} matrix to all of the available data: if there
are $n_{ij}$ transitions observed from state $i$ to state $j$, then the
corresponding element of the fitted transition matrix is
$p_{ij}=n_{ij}/\sum_{k=0}^m n_{ik}$.
A small refinement is to observe that for
\textcolor{black}{the Markov chain $(Z_t)$} the stationary probability
distribution is already known: the stationary probability $\pi_i$ of
state $i$ should be $p_i - p_{i-1}$ ($i=1,\ldots,m$).
In the special case of equally spaced quantiles ($p_i=i/m$), the stationary
distribution is uniform ($\pi_i=1/m$ for each $i$); equivalently, the 
transition matrix is doubly stochastic.
A transition matrix constrained to satisfy this requirement can be estimated
{\it e.g.} by the Sinkhorn-Knopp algorithm (\cite{Sinkhorn,SinkhornKnopp}).
Such a model is described in detail in \cite{Rayner}.

More ambitiously, we can consider that the serial dependence structure of the original process
may also exhibit variation throughout the periodic cycle.
To accommodate such variation, the transition probabilities $p_{ij}(t)$ should be
allowed to be time-inhomogeneous: periodic functions of $t$, with the same period $T$
as the original process. That is,
\begin{equation}
\label{eq:periodic_pij}
 p_{ij}(t) = \sum_{\ell=0}^q \gamma_{ij\ell} b_\ell(t)
\end{equation}
where $b_0(t),\ldots,b_q(t)$ are periodic basis functions and $\gamma_{ij\ell}$ are fitted coefficients.

This is a separate and different kind of periodicity from that contemplated in Section \ref{s:QFR}:
while (\ref{eq:QR}) describes a periodic variation in the absolute meaning of the Markov states,
(\ref{eq:periodic_pij}) describes a periodic variation in the transition behaviour between states.
Suppose, for example, that the original process $(X_t)$ represents daily rainfall; then the
highest Markov state $m$ represents a rainy day. For any time of year $t$, the function
$q_{p_{m-1}}(t)$ gives information on the likely amount of rain falling on one rainy day,
while $p_{mm}(t)$ gives information on the likelihood of a sequence of consecutive rainy days.
\textcolor{black}{(In the time-homogeneous case, the length of such a
sequence would have a geometric distribution with parameter $p_{mm}(t)$.)}

The periodic function space used in (\ref{eq:periodic_pij}) can be the same as that used in
(\ref{eq:QR}), or different. For the examples in the present paper, we will again use the
Fourier basis.

\medskip
\textcolor{black}{
{\bf Example.} Suppose that a two-state process with
period 12 ({\it e.g.} monthly observations and annual period) is to be
modelled as a stationary time-inhomogeneous Markov chain with transition
matrix at time $t$ given by
$$ P(t) = (p_{ij}(t))_{i,j=1}^2 =
\begin{pmatrix} \gamma_{110}, \gamma_{120} \\ \gamma_{210}, \gamma_{220}
  \end{pmatrix}
+
\begin{pmatrix} \gamma_{111}, \gamma_{121} \\ \gamma_{211}, \gamma_{221}
  \end{pmatrix} \cos(\omega t)
+
\begin{pmatrix} \gamma_{112}, \gamma_{122} \\ \gamma_{212}, \gamma_{222}
  \end{pmatrix} \sin(\omega t) ,
  $$
where $\omega=\pi/6$.
Since the entries of $P(t)$ are probabilities,
$$ \gamma_{ij0} + \gamma_{ij1}\cos(\omega t) + \gamma_{ij2}\sin(\omega t)
\in[0,1] \qquad\hbox{ for $i,j=1,2$ and all $t$,}
  $$
which can be expressed via the constraints
\begin{equation}
\label{eq:example_ineq_constraints}
 \left(\gamma_{ij1}^2 + \gamma_{ij2}^2\right)^{1/2}
   \leq \hbox{min}(\gamma_{ij0},1-\gamma_{ij0})
\qquad\hbox{ for $i,j=1,2$.}
\end{equation}
The rows of $P(t)$ sum to 1 at all times $t$, giving the linear constraints
\begin{eqnarray}
\label{eq:example_row_constraints}
&& \gamma_{110} + \gamma_{120} = 1  \quad\hbox{ and }\quad
   \gamma_{111} + \gamma_{121} = 0, \quad
   \gamma_{112} + \gamma_{122} = 0\\
&& \gamma_{210} + \gamma_{220} = 1  \quad\hbox{ and }\quad
   \gamma_{211} + \gamma_{221} = 0, \quad
   \gamma_{212} + \gamma_{222} = 0.
\end{eqnarray}
If we require that $P(t)$ be doubly stochastic (so that the two states
are equally likely at any given time $t$) then its columns must also
sum to 1, giving
\begin{eqnarray}
&& \gamma_{110} + \gamma_{210} = 1  \quad\hbox{ and }\quad
   \gamma_{111} + \gamma_{211} = 0, \quad
   \gamma_{112} + \gamma_{212} = 0\\
&& \gamma_{120} + \gamma_{220} = 1  \quad\hbox{ and }\quad
   \gamma_{121} + \gamma_{221} = 0, \quad
   \gamma_{122} + \gamma_{222} = 0.
\label{eq:example_col_constraints}
\end{eqnarray}
Note that the linear equality constraints
(\ref{eq:example_row_constraints}--\ref{eq:example_col_constraints})
are not all linearly independent.
}  

\medskip
{\bf Estimation of coefficients.}
Estimation of the $\gamma_{ij\ell}$ can be done by a maximum
likelihood procedure.
\textcolor{black}{Suppose that for each $k=1,\ldots,n$ we have observed a transition at time $t_k$ from a state $i_k$ to a state $j_k$.
(The notation $t_k$ for the time of the $k$th observation allows for the
possibility of missing or non-consecutive data; in the absence of
such complications $t_k=k$.) The likelihood (probability) of this
collection of observations}
is $L=\prod_{k=1}^n p_{i_k j_k}(t_k)$.
\textcolor{black}{The coefficients $\gamma_{ij\ell}$ in the model should thus be}
chosen to minimize
$$ -\log L = -\sum_{k=1}^n \log\left(\sum_{\ell=0}^q \gamma_{i_k j_k \ell} b_\ell(t_k)\right)
  $$
subject to the constraints
\begin{eqnarray*}
0 \leq \sum_{\ell=0}^q \gamma_{ij\ell} b_\ell(t) \leq 1 & \qquad\hbox{ for all $i$, $j$, and $t$}\\
\sum_{j=1}^m \sum_{\ell=0}^q \gamma_{ij\ell} b_\ell(t) = 1 & \qquad\hbox{ for all $i$ and $t$}\\
\sum_{i=1}^m \pi_i \sum_{\ell=0}^q \gamma_{ij\ell} b_\ell(t) = \pi_j & \qquad\hbox{ for all $j$ and $t$}
\end{eqnarray*}
where $\pi_i$ is the known stationary probability of state $i$ noted above.
This is a nonlinear convex optimization problem
\textcolor{black}{which can be solved by any of a number of standard methods}.
Note the requirement that the constraints hold
for all times $t$, rather than for only those times $t_k$ at which transitions have been observed.
In the case of the Fourier basis (\ref{eq:Fourier_basis}), the equality constraints
reduce to
\begin{eqnarray*}
\sum_{j=1}^m \gamma_{ij0} = 1 \quad\hbox{ for all $i$}\qquad&\hbox{ and }&\qquad
   \sum_{j=1}^m \gamma_{ij\ell} = 0 \quad\hbox{ for all $i$ and $\ell>0$}\\
\sum_{i=1}^m \pi_i \gamma_{ij0} = \pi_j \quad\hbox{ for all $j$}\qquad&\hbox{ and }&\qquad
   \sum_{i=1}^m \pi_i \gamma_{ij\ell} = 0 \quad\hbox{ for all $j$ and $\ell>0$}.
\end{eqnarray*}

\medskip
\textcolor{black}{
{\bf Example continued.}
Suppose that the following sequence of states is observed.
$$ 1 2 2 2 2 2 2 1 2 1 2 1 2 2 2 2 2 2 1 2 1 2 1 2 2 1 1 1 1 1 1
  $$
That is, the first observed transition is from state $i_1=1$ to state
$j_1=2$ at time $t_1=1$; the likelihood of this is its probability
$$ p_{12}(1)
= \gamma_{i_1 j_1 0} + \gamma_{i_1 j_1 1}\cos(\omega t_1)
  + \gamma_{i_1 j_1 2}\sin(\omega t_1)
= \gamma_{120} + \gamma_{121}\cos(\omega) + \gamma_{122}\sin(\omega)
  $$
The last observed transition is from state $i_{30}=1$ to state
$j_{30}=1$ at time $t_{30}=30$; its likelihood is $p_{11}(30)$.
The likelihood of the complete sequence is the product
$$ L = p_{12}(1) \cdot p_{22}(2)\cdots p_{11}(30)
  $$
and the objective to be minimized is
$$ -\log L = -\sum_{k=1}^{30}
  \log\left(\gamma_{i_k j_k 0}
          + \gamma_{i_k j_k 1}\cos(\omega t_k)
          + \gamma_{i_k j_k 2}\sin(\omega t_k)\right)
  $$
subject to the constraints (\ref{eq:example_ineq_constraints}--\ref{eq:example_col_constraints}).
}  

\section{A Markov decision process}
Equipped with a Markov process that describes an underlying natural phenomenon, {\it e.g.} rainfall,
we can proceed to the decision-making part of the model. Markov decision processes span a wide range of models,
that include finite or infinite state as well as discrete or continuous processes.
For the purposes of this paper, we focus on finite discrete MDPs. {\textcolor{black}{We will consider an infinite time horizon MDP where we obtain policies that minimize the long-run average per-cycle cost. As with any Markov decision process, our MDP has the following four main components:
\begin{enumerate}
    \item A set $S$ of states, with members $i\in S$. 
    \item A set $K$ of actions with $k \in K$. 
    \item Costs or rewards associated with taking an action $k$ in state $i$, at time $t$. 
    \item Probability of transitioning from state $i$ to state $j$, when action $k$ is undertaken, at time $t$, denoted by $P_{ijt}(k)$.
\end{enumerate}
}}
\textcolor{black}{We now describe these components in more detail. The state space $S$ for the MDP consists of the Cartesian product of states of an underlying natural phenomenon and system states that are partly, or wholly driven by allowable actions.
For the reservoir management application in the next section,
the natural phenomenon is inflow into the reservoir. We model this as a Markov
process and quantiles of the weekly natural inflow into the Waitaki River
system define the states of this Markov chain. 
For the MDP, a state $i \in S$ is a tuple $i = (q, l)$ where $q$ is the quantile of natural inflow, and $l$ is the (discretized) storage level of the reservoir (note that the reservoir level is a byproduct of natural inflow {\em and} the utilization of the reservoir to produce energy).
In other, similar, energy-related problems, MDP states may include the levels of
charge stored in batteries, or the internal states of energy demand models.}
We will assume that at any time $t$ and in any given state $i$,
we can take a range of actions (or decisions) $k \in \{0,1, \cdots, |K|-1\}$.
\textcolor{black}{In our reservoir management application, our actions amount
to the dispatch of non-hydropower resources. Non-hydro dispatch amounts are
discretized in 100MW blocks, so the action set is denoted by $k \in \{0, 100,
200, \cdots, 800\}$. Any shortfall from electricity demand is then met through electricity production from the reservoir.}

The effect of taking an action is twofold: firstly, it incurs a cost;
secondly, it determines the probability distribution of the state at time
$t+1$.
\textcolor{black}{The cost of taking action $k$ when the system is in state $i$
at time $t$ is denoted by $c_{ikt}$. It is also possible to assign a cost
to merely visiting a state $i$, by including that cost in $c_{ikt}$ for all
actions $k$ possible for that state.
For the reservoir management example, a cost associated with curtailment of
demand (value of lost load) is incurred in any state where demand exceeds
available supply, regardless of the action taken.}

\textcolor{black}{The probability of transition to state $j$
when the system is in state $i$ at time $t$ and action $k$
is denoted $P_{ijt}(k)$. These transition probabilities can be derived
from the Markov model of the underlying natural phenomenon,
as this is the only source of randomness in the system.}

\textcolor{black}{The reader should note that as in any other MDP, the states and actions are tailored to the underlying application. We specify these in detail for our case studies below.}

Markov decision processes select policies that lead to a desired objective,
such as minimizing long-run expected average cost per unit time.
A randomized policy will present the decision maker with a probability distribution
$(d_{i1t}, d_{i2t}, \cdots, d_{iKt})$ for the action $k$ to be chosen when in state $i$ at time $t$. 
Given such a policy, the system can operate as a cyclostationary process in the sense given
in Section \ref{s:QFR}; this represents its long-run behaviour under the specified policy.
We can find the policies minimizing the average (over the periodic cycle) expected cost
by forming the following linear program. {\textcolor{black}{This is the model that we follow for our case studies below.}} Note that choosing an action affects the
probability of reaching other states as well as the immediate cost accrued, both of which will have
an impact on the steady state distribution of cost, hence the expectation of cost.

\begin{equation}
\label{eq:MDPLP}
\begin{array}{lll}
\mbox{Minimize} & \sum_t \sum_{i=0}^M\sum_{k=0}^K c_{ikt}y_{ikt} & \\
\mbox{s/t}      & \sum_{i=0}^M\sum_{k=0}^K y_{ikt} = 1 & \forall t \\
                 & \sum_{k=0}^K y_{j,k,t+1} - \sum_{i=0}^M\sum_{k=0}^K y_{ikt} P_{ijt}(k) = 0 & \forall t, \forall j\\
                 & y_{ikt} \geq 0 & \forall i, \forall k, \forall t.
\end{array}
\end{equation}

This linear program is an extended version of those given in \cite{HillierLieberman, White},
where we have introduced periodic time inhomogeneity, hence the dependence on $t$. 
As previously stated, the parameter $P_{ijt}(k)$ represents transition probabilities.
The parameter $C_{ikt}$ is the cost incurred by being in state $i$ at time $t$ and taking action $k$.
Decision variable $y_{ikt}$ represents the cyclostationary probability that the system is in
state $i$ and action $k$ is undertaken at time $t$.
Once the optimal $y_{ikt}$ variables are obtained, we can determine $d_{ikt}$ using
$d_{ikt} = \frac{y_{ikt}}{\sum_{k=0}^K y_{ikt}}$.

\section{Application: reservoir management}
\label{s:hydro}

In this section, we apply our technique to a hydro-power reservoir management
problem of the kind described in \cite{WangAdams}.

\textcolor{black}{For our example problem, we require an electric power system
to meet a constant demand for 1400\,MW of electric power at the lowest
possible average cost. The supply resources available comprise 800\,MW of
fossil-fueled power generation capacity and a hydropower generation system
modelled loosely on the Waitaki River system in New Zealand.
The hydropower generation capacity is taken to be 1500\,MW.
Assumed costs are \$50 per megawatt-hour (MWh) for non-hydro generation
and \$1000/MWh for unserved demand. The hydropower has no cost, but is
limited by available water inflows.}

\textcolor{black}{The state space of this problem comprises the internal state
of a four-state Markov chain model of hydro inflows (see below) together with
the energy storage state of the hydro reservoir.
Although the real Waitaki River system comprises multiple hydrological
catchments and reservoirs, for the purposes of this paper it may be treated
as a single equivalent-energy reservoir. The energy storage capacity
of the real system is approximately 2500 gigawatt-hours (GWh)
\cite{EA_HMD}, but we will here consider an illustrative problem in which the
reservoir capacity is only 840\,GWh. (Downsizing the reservoir is done
principally to make the problem more interesting, by increasing the
probabilities of shortages and overflows.) The Markov decision process
requires discrete states, which we create by discretizing the stored energy
into blocks of 16.8\,GWh. Since the time step is one week, this block size can
be conveniently expressed in power units as 100\,MW. The reservoir capacity is
50 blocks, and so there are 51 possible storage states.
With the four possible inflow states, the size of the state space is 204.}

\textcolor{black}{We seek a cyclostationary solution with annual period,
consisting of $52\times 7$-day time steps, for which the expected cost per
annum of non-hydro power and unserved demand is minimized.
At each time step, the available actions consist of dispatching some amount
of non-hydro generation; since this must be delivered in 100\,MW blocks
from a total capacity of 800\,MW, there are 9 different actions available.
Thus, the variables $y_{ikt}$ in (\ref{eq:MDPLP}) number
$204\times9\times52=95472$.}

\begin{figure}
\begin{center}
\centerline{
\includegraphics[width=0.5\linewidth]{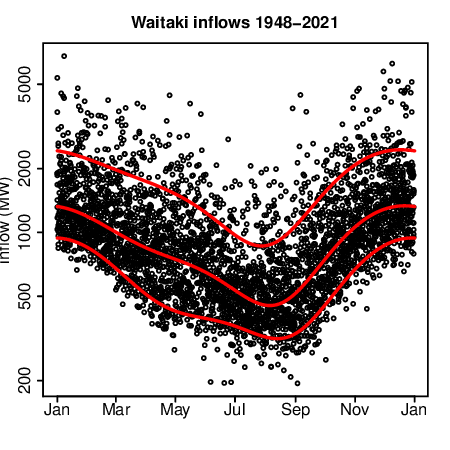}
\hfil
\includegraphics[width=0.5\linewidth]{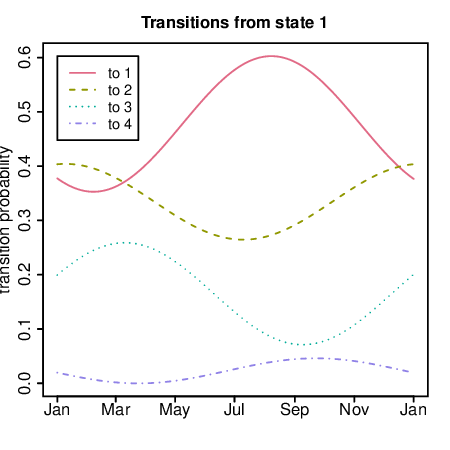}
}
\caption{The Waitaki River model. Left panel: weekly inflow data, with fitted 10th, 50th, and 90th percentile models \textcolor{black}{partitioning the
data into four inflow states}.
Right panel: \textcolor{black}{Transition probabilities from the first (lowest-flow) Markov state, fitted acording to the methodology in Section \ref{s:Markov_model}}.}
\label{fig:hydro_quantile_model}
\end{center}
\end{figure}

\textcolor{black}{The Markov inflow model is derived via the quantile Fourier regression methodology laid out in Section \ref{s:QFR}.
The underlying data comprise a univariate time series of weekly natural
inflows to the Waitaki River system; see \cite{Pritchard_hydro}
for more details on this data set.
For convenience, we have expressed the inflows in power units (megawatts).}
Figure \ref{fig:hydro_quantile_model} (left panel) shows the data set together with quantile regression models
for the annual variation of the 10th, 50th, and 90th percentiles. These models were constructed with the
Fourier basis (\ref{eq:Fourier_basis}) with $r=2$; that is, the fitted curves are second-order trigonometric polynomials.

These three quantile models were then used to model the serial dependence structure as
a four-state Markov chain as described in Section \ref{s:Markov_model}.
All sixteen transition probabilities were permitted to vary annually as simple sinusoids,
by using the Fourier basis (\ref{eq:Fourier_basis}) with $r=1$.
The fitted probabilities for transitions from state 1 (lowest flow) to other states are shown
in Figure \ref{fig:hydro_quantile_model} (right panel).

The Markov inflow process must now be further developed to deliver energy inflows in blocks
of the same size as used for the reservoir storage.
Each state $i$ of the original process discussed in Section \ref{s:Markov_model} corresponds only to
a time-varying interval $[q_{p_{i-1}}(t), q_{p_i}(t)]$ containing the amount of energy inflow.
Our chosen block size is small enough that this interval usually contains least two
different multiples of the block size. So, we can take the inflow amount to be
one of these multiples, chosen at random independently of the underlying Markov process.
Note that this approach preserves the continuous flavour of the original inflow model:
even though the inflow amounts are now discrete, their probabilities can be thought of as
varying continuously with time $t$.

\begin{figure}
\begin{center}
\centerline{
\includegraphics[angle=-90,origin=c,width=\linewidth]{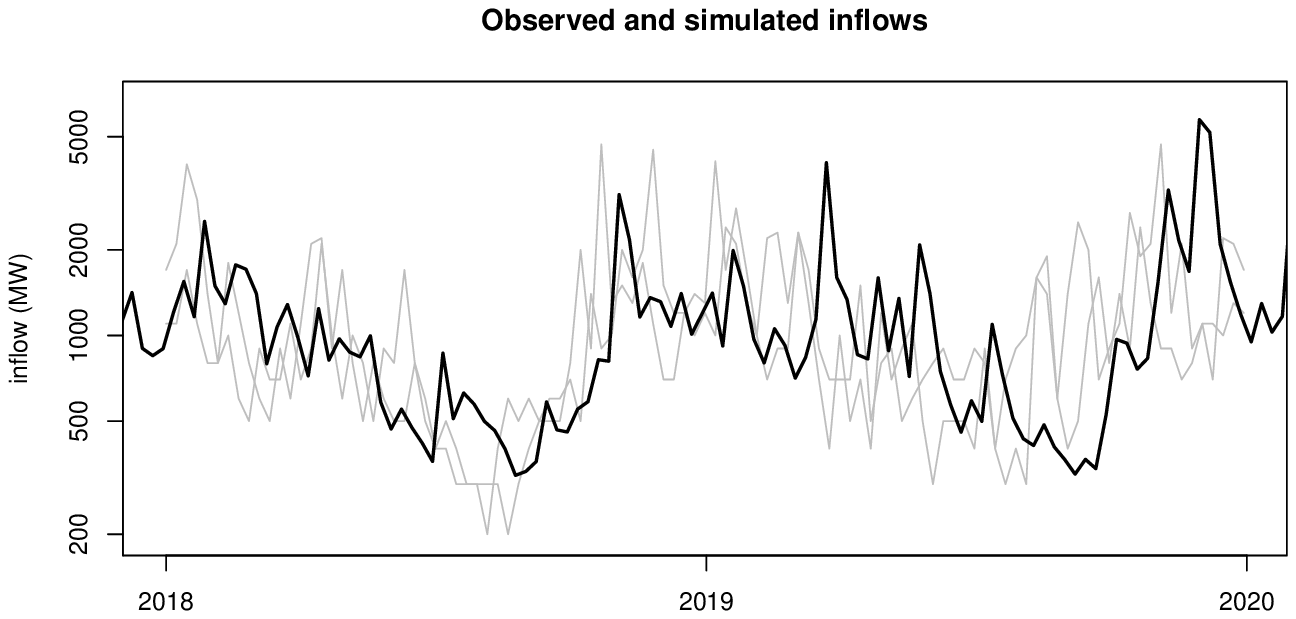}
}
\vspace{-100pt}
\caption{\textcolor{black}{Two simulated sample paths of Waitaki River inflows
in 2018 and 2019 (in grey). The darker line shows actual inflows over
the same period.}}
\label{fig:inflow_simulation}
\end{center}
\end{figure}

\begin{figure}
\begin{center}
\centerline{
\includegraphics[width=0.5\linewidth]{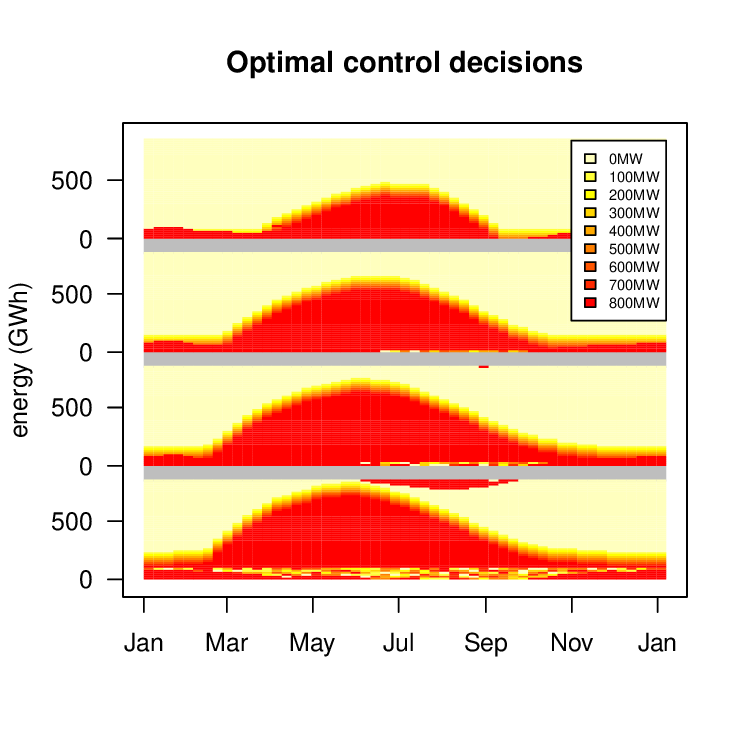}
\hfil
\includegraphics[width=0.5\linewidth]{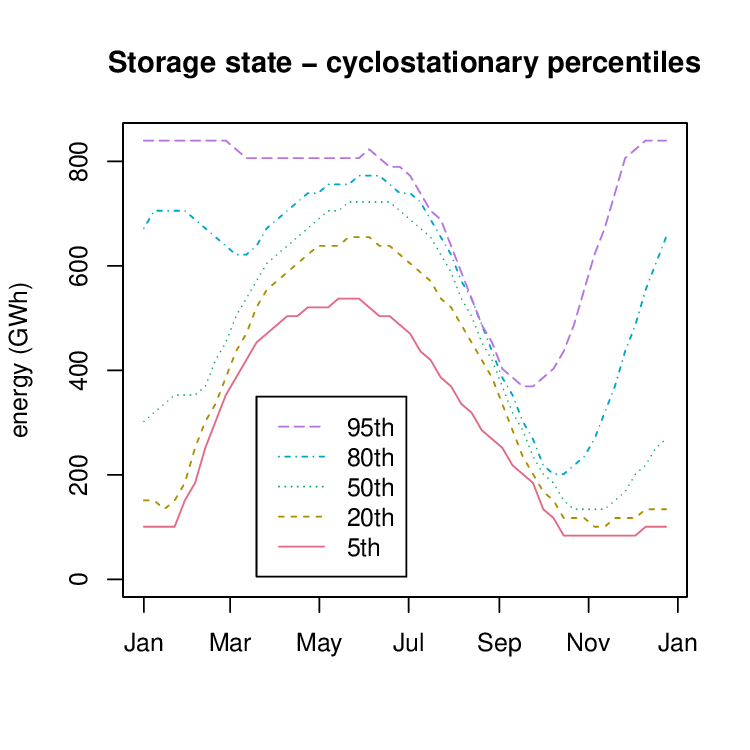}
}
\caption{The Waitaki River MDP solution. Left panel: the optimal decisions (non-hydro generation quantities) \textcolor{black}{as a function of time
(annual cycle), stored energy (0-840\,GWh), and inflow state (four states)}.
Right panel: the \textcolor{black}{probability distribution of stored energy
as a function of time for the cyclostationary process}.}
\label{fig:hydro_solution}
\end{center}
\end{figure}

We determine the distribution of the inflow, conditional on state $i$ and time $t$, in the
following way. For states other than the highest and lowest, the original inflow data
corresponding to state $i$ are normalized by the transformation
$$ x \mapsto \frac{x - q_{p_{i-1}}(t_x)}{q_{p_i}(t_x) - q_{p_{i-1}}(t_x)}
  $$
(where $t_x$ is the time at which inflow $x$ was observed).
\textcolor{black}{This gives an empirical probability
distribution on $[0,1]$, estimating the conditional distribution of
normalized inflow given the state $i$.}
The multiples of the block size permissible at time $t$ are then
transformed in the same way, and assigned probabilities to match the empirical distribution
as closely as possible. For the lowest and highest states, the corresponding intervals
are semi-infinite: $(-\infty, q_{p_1}(t)]$ and $[q_{p_{m-1}}(t),\infty)$.
In these cases, the normalization transformation used is a simple division by the finite endpoint.

\textcolor{black}{Figure \ref{fig:inflow_simulation} illustrates
inflow sample paths generated by this time-inhomogeneous Markov model.}

Figure \ref{fig:hydro_solution} illustrates the solution of the MDP.
Relatively abundant summer inflows create a wide range of possible storage states at the beginning of the
year. But by the beginning of May (late autumn), the optimal control has used non-hydro generation to
ensure that the storage reservoir is fairly full, so that the stored energy can be drawn upon during the
low-inflow winter period.

\section{Application: thermal backup of offshore wind}
\label{s:offshore_wind}

With the increasing penetration of renewable sources of electricity generation into the world's electricity systems,
it is imperative to have sufficient backup (or firming) for intermittent renewables and to operate it efficiently
to avoid costly curtailing of demand for electricity.
In the United States there are a number of mandates to procure electricity generation from offshore wind;
for instance, Massachusetts has set a target of 3.2 gigawatts (GW) of offshore wind capacity by 2035 \cite{Massachusetts}.

Our first example concerns the operation of a system of backups to cover any shortfall in power generated
from offshore wind resources. The model's underlying univariate stochastic process represents the New England
regional demand for electric power \textcolor{black}{minus offshore wind power generation. The figures for wind power generation are computed from wind speed data (below we specify this calibration process).}

\begin{figure}
\begin{center}
\centerline{
\includegraphics[width=0.5\linewidth]{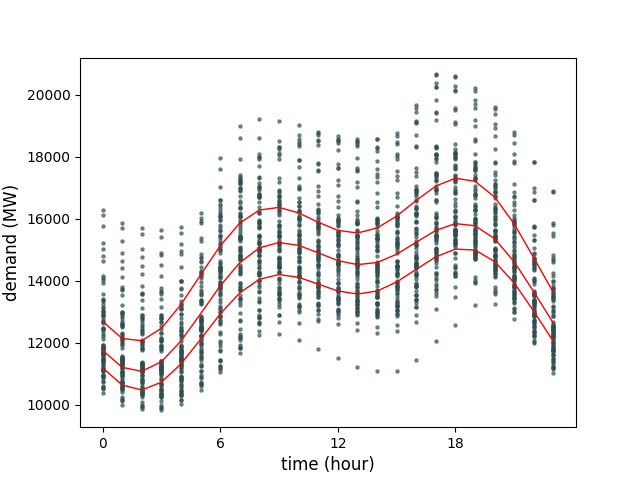}
\hfil
\includegraphics[width=0.5\linewidth]{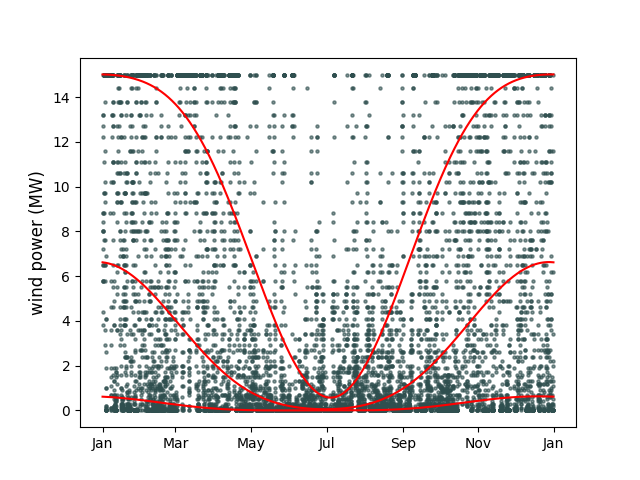}
}
\caption{Demand and wind power processes, with 25th, 50th, and 75th percentiles estimated by quantile Fourier regression with $r=2$
and superimposed on 2018 data. Left panel: diurnal variation in New England winter demand.
Right panel: annual variation in generation by one IEA 15\,MW offshore wind turbine.}
\label{fig:demand_and_wind}
\end{center}
\end{figure}

We use hourly electricity demand data for the ISO-New England grid from 2006--2020
as reported to the Federal Energy Regulatory Commission (FERC) \cite{FERC}.
Figure \ref{fig:demand_and_wind} (left panel) illustrates the diurnal variation,
showing a double-peaked pattern typical of electricity consumption
\textcolor{black}{in the winter months (Dec, Jan, Feb), but not of summer.} 
This pattern indicates that quantile modelling with annual period must also capture daily periodicity
to at least the second harmonic and that the shape of the diurnal variation should be modulated by the time of year.
That is, the basis functions in (\ref{eq:QR}) should include {\it e.g.} 
$\cos(\omega_1 t)\cos(2\omega_2 t)$ and similar functions,
where $\omega_1=2\pi/\hbox{year}$ and $\omega_2=2\pi/\hbox{day}$.
\textcolor{black}{A full enumeration of this basis is given in the Appendix
to this paper.}

For wind power, the underlying data are hourly wind speeds recorded by the
National Oceanic and Atmospheric Administration (NOAA) National Data Buoy Center \cite{NOAA}.
As a reference station for evaluating offshore wind power generation, we used Buoy 44025, in the New York Bight area.
We chose the IEA 15\,MW reference turbine \cite{IEA15} as our power generation mechanism, using the power curve 
and other specifications (hub height, etc.) for this turbine to make the wind-speed to power conversion.
Details of the power curve modeling can be found in \cite{Betz, Kirchhoff}.
Figure \ref{fig:demand_and_wind} (right panel) illustrates the use of quantile Fourier regressions for wind power
generated at this scale. It is clear that there are seasonal patterns in wind power,
with a significant drop in wind generation in the summer.
To calibrate the contribution of offshore wind, we were guided by the ISO-New England target for wind generation for 2030.
Based on \cite{NEGoal,isoNE}, we assume 21687\,MW of offshore windpower capacity, corresponding to
1446 IEA-15\,MW turbines in the region.

\begin{figure} [ht!]
\centering
\begin{minipage}{.5\textwidth}
  \centering
  \includegraphics[width=1\textwidth]{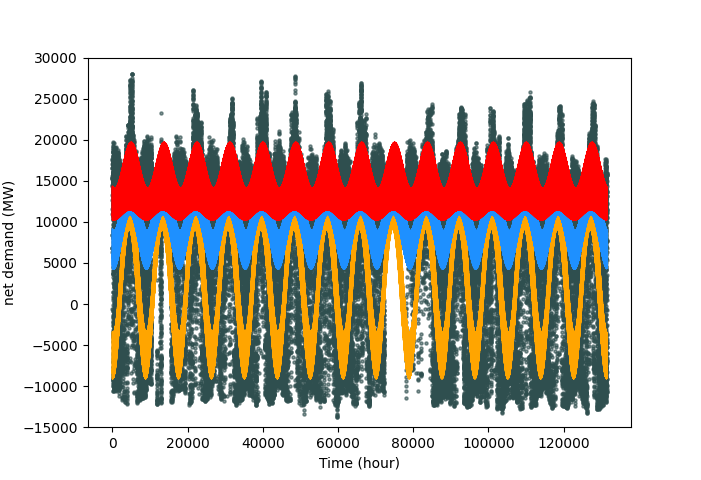}
\end{minipage}%
\begin{minipage}{.5\textwidth}
  \centering
  \includegraphics[width=1\textwidth]{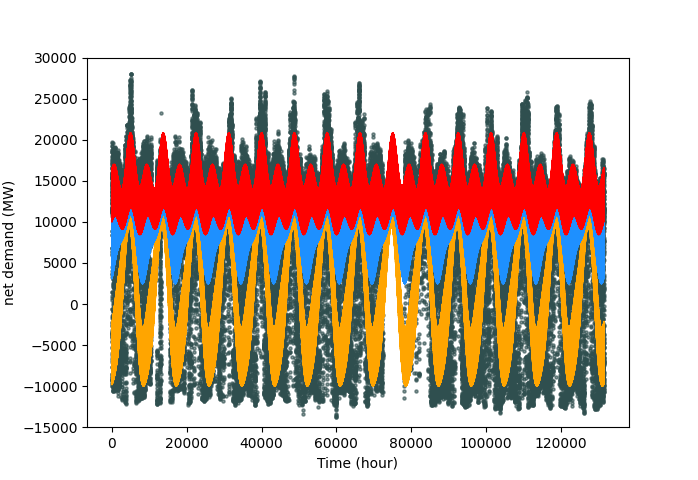}
\end{minipage}
\caption{Quantile Fourier regression with both annual and daily periodicity, fitted to 15 years of net demand data.
Left panel: $r=1$. Right panel: $r=2$.}
\label{fig:net_demand_quantiles_15_years}
\end{figure}

\begin{figure} [ht!]
\centering
\begin{minipage}{.5\textwidth}
  \centering
  \includegraphics[width=1\textwidth]{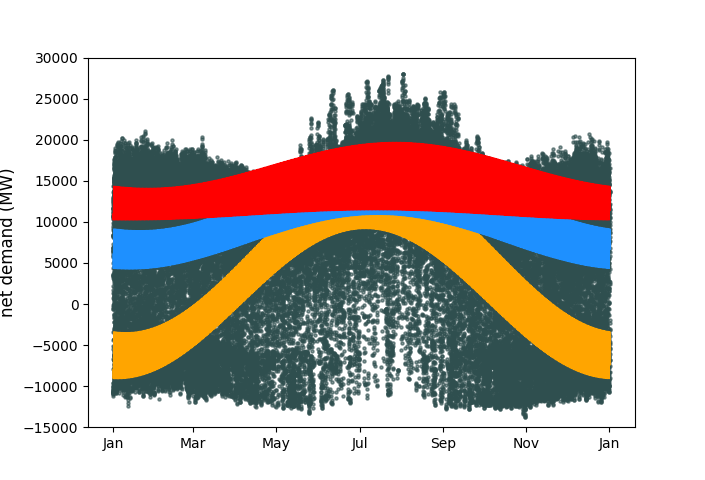}
\end{minipage}%
\begin{minipage}{.5\textwidth}
  \centering
  \includegraphics[width=1\textwidth]{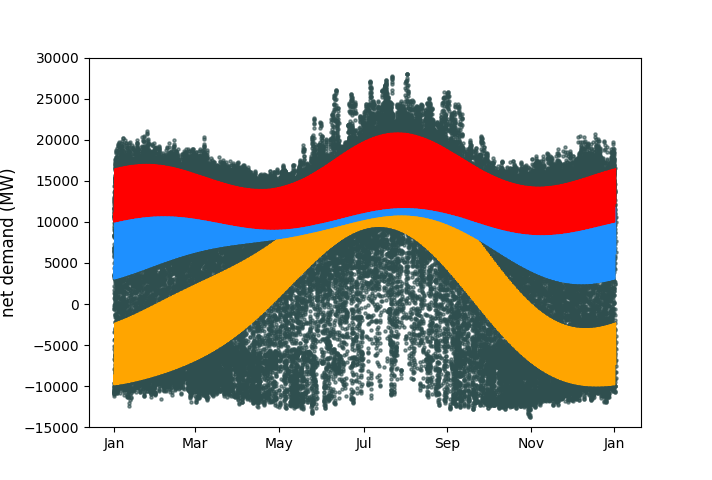}
\end{minipage}
\caption{Quantile Fourier regression with both annual and daily periodicity, fitted to 15 years of net demand data
(phase-folded plots). Left panel: $r=1$. Right panel: $r=2$.}
\label{fig:net_demand_quantiles_phase_folded}
\end{figure}

\begin{figure}[ht!]
\begin{center}
\centerline{
\includegraphics[width=0.5\linewidth]{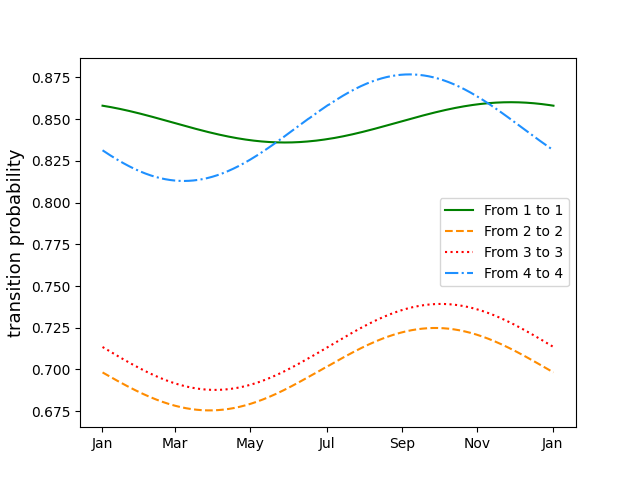}
\hfil
\includegraphics[width=0.5\linewidth]{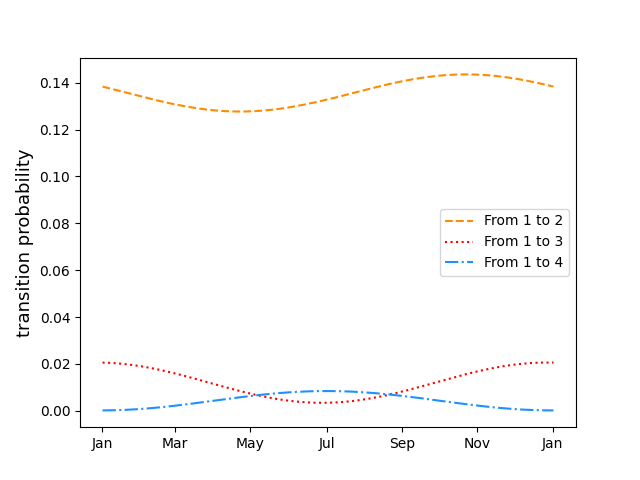}
}
\caption{Fitted transition probabilities for the net demand process. Left panel: transitions from each state to itself.
Right panel: transitions from the first (lowest) net-demand state to other states.}
\label{fig:net_demand_transition_probs}
\end{center}
\end{figure}

Since our model is concerned only with the difference between power demand and wind power supply,
we combine the two into a single univariate time series of hourly net demand
before proceeding to the modelling stage.
Figures \ref{fig:net_demand_quantiles_15_years} and \ref{fig:net_demand_quantiles_phase_folded}
juxtapose quantile Fourier regression fits to the net demand series.
The improvement obtained by including the second harmonic ($r=2$) is apparent.
\textcolor{black} {Visual inspection demonstrates a better fit for $r=2$ as the peaks and off-peaks are better represented. We also see this in the pseudo $R^2$ for the quantile regression, which is 0.19, 0.10, and 0.12 for quantiles of 0.25, 0.5, and 0.75 respectively of $r=1$, while it is 0.48, 0.48, and 0.49 for quantiles of 0.25, 0.5, and 0.75 respectively of $r=2$.
While we can increase $r$ ad infenitum, as $n$ increases, the number of terms in the Fourier regression fit increases quadratically. Given this rapid non-linear increase, we would prefer smaller $r$s sufficient for rendering near optimal operational plans. There is also a natural intuition behind choosing $r=2$, as it will allow adapting our fits to day/night partition of a 24 hour period (Figure \ref{fig:demand_and_wind} left panel), similarly annual fits to wind power would naturally benefit from distinguishing summer and winter (Figure \ref{fig:demand_and_wind} right panel).
}

The quantile models were then used to model the serial dependence structure as
a four-state Markov chain {\textcolor{black}{as described in Section \ref{s:Markov_model} of this paper.}
All sixteen transition probabilities were permitted to vary annually as simple sinusoids,
by using the Fourier basis (\ref{eq:Fourier_basis}) with $r=1$.
The resulting transition probabilities are illustrated in Figure \ref{fig:net_demand_transition_probs}.

Building on the net demand model, we constructed a Markov decision process containing a stylized
representation of fourteen combined-cycle gas turbine (CCGT) powerplants with a total of 28000\,MW 
generation capacity.
For this stylized example, in each state, we allow increased generation of 2000\,MW as ``ramp up" action
(provided thermal generation level in the current state is 26000\,MW or less), decrease of generation by 2000\,MW as 
``ramp down" action (provided current thermal generation is at least 2000\,MW), or staying at the same thermal 
generation level.
Equivalently, we have an enormous thermal power plant that can ramp up or down by 2000\,MW in each period.
The model could be made closer to reality by allowing multiple CCGTs to ramp in each state;
this would increase the number of allowable actions in each period significantly. 

This example allows us to illustrate the effect of time inhomogeneity in a clear and concise manner.
Figure \ref{fig:ccgt_decision_plot_non_seasonal} illustrates the optimal operation plans that result
from a Markov decision process model in which the decisions may vary by time of day, but not by time of year.
In contrast Figure \ref{fig:ccgt_decision_plot_seasonal} demonstrates optimal decisions for the full
model, in which the operation plan may be different on each day of the year. 
In these figures, red indicates ``ramp up" is the recommended action for the state,
yellow indicates no change, and green is indicative of ``ramp down".
Note the significant differences between winter (upper panel) and summer (lower panel) operations.
In particular, the thermal back up operation level rises to cover the afternoon peak demand during the summer,
which could be attributed to air conditioning use in hot July afternoons,
together with the decrease in wind power generation in the summer.
Both operational plans laid out in Figure \ref{fig:ccgt_decision_plot_seasonal} differ from the all-year-round 
plan in Figure \ref{fig:ccgt_decision_plot_non_seasonal}.

\begin{figure}[ht!]
\centering
\includegraphics[width=1\textwidth]{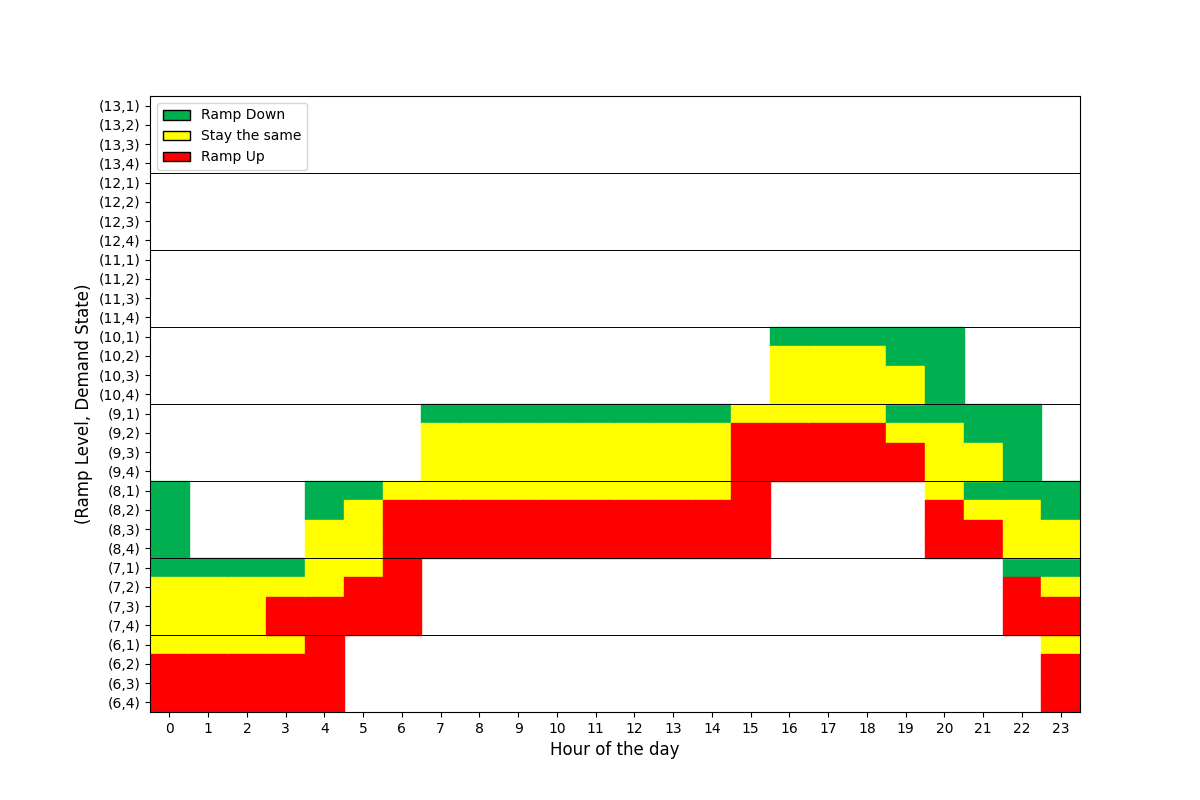}
\includegraphics[width=1\textwidth]{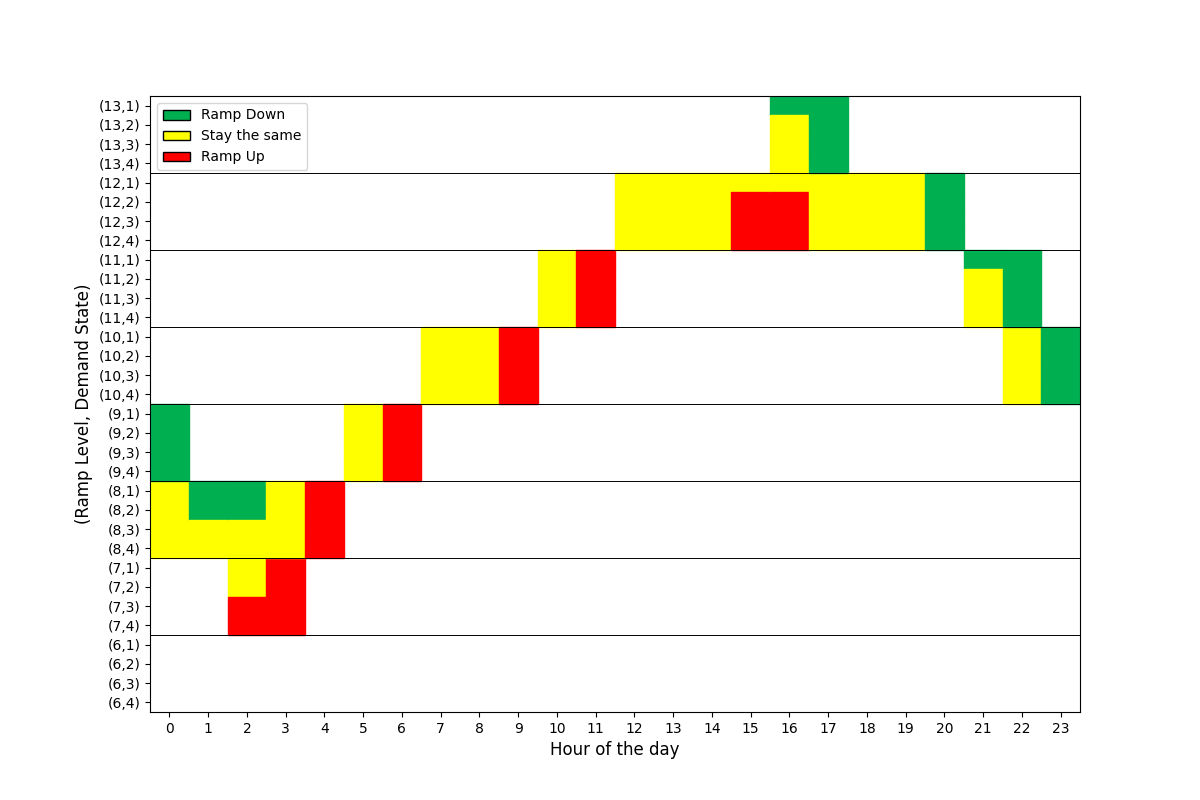}
\caption{The offshore wind integration MDP solution. Upper panel: optimal decisions on January 1.
Lower panel: optimal decisions of July 1. Here the colors indicate the action to be taken in each 
state at each time: red for ramp up, green for ramp down and yellow for status quo.}
\label{fig:ccgt_decision_plot_seasonal}
\end{figure}

\begin{figure} [ht!]
\centering
  \includegraphics[width=1\textwidth]{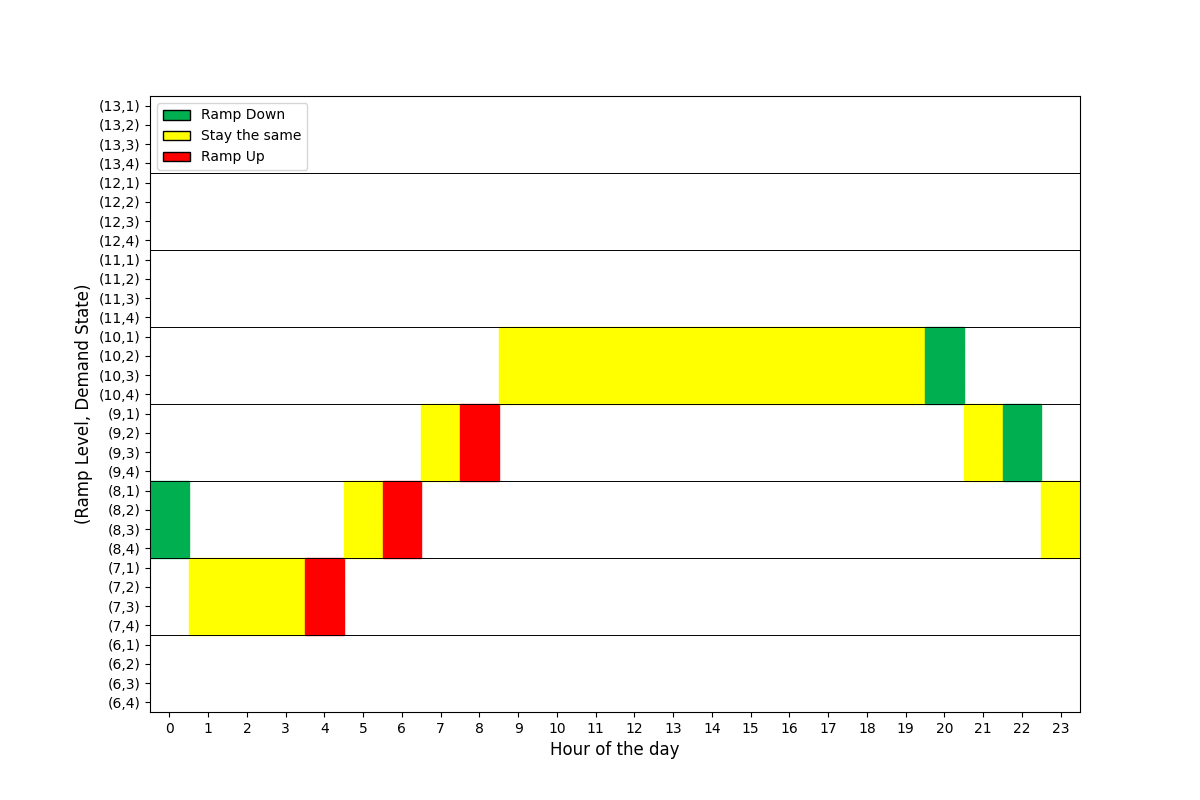}
\caption{The offshore wind integration MDP solution: optimal decisions when the decision rule is required to be the same on every day of the year.}
\label{fig:ccgt_decision_plot_non_seasonal}
\end{figure}

\begin{table}[ht!]
\centering
\caption{Comparison of the results of a simulation between time-dependent plan and fixed plan over 13 years of data. (Extra is the extra demand that is still unmet after CCGTs production)}
\label{tab:cost_comparison}
\begin{tabular}{|c|c|c|}
\hline
\textbf{Approach} & \textbf{Total Extra} & \textbf{Total Cost} \\
\hline
Time-dependent action-plan & 13,776 & 186,207 \\
\hline
Fixed plan for every day of the year & 33,534 & 245,736 \\
\hline
\end{tabular}
\end{table}

\section{Acknowledgements and compliance with ethical standards}
The authors would like to acknowledge the National Science Foundation (GCR award 2020888), The Sloan Foundation (award number 2023-19608) and ISO New England for their generous support of our research. 
This article does not contain any studies with human participants or animals performed by any of the authors.

\newpage
\section{Appendix}
Periodic time-varying quantile functions with both annual
and diurnal frequencies, as used in Section \ref{s:offshore_wind}},
have the following form when $r=1$.

\begin{equation}
\label{eq:FR1}
\begin{aligned}
q_{p}(t) &= \mu + A_1\cos{\omega_{day}t} + B_1 \sin{\omega_{day}t} \\
&\quad \quad + A'_1 \cos{\omega_{year}t} + B'_1 \sin{\omega_{year}t} \\
&\quad \quad + C_1 \cos{(\omega_{year}t})\cos{(\omega_{day}t}) \\
&\quad \quad + D_1 \cos{(\omega_{year}t})\sin{(\omega_{day}t}) \\
&\quad \quad + E_1 \sin{(\omega_{year}t})\cos{(\omega_{day}t}) \\
&\quad \quad + F_1 \sin{(\omega_{year}t})\sin{(\omega_{day}t})
\\
\end{aligned}
\end{equation}

When $r=2$,

\begin{equation}
\label{eq:FR2}
\begin{aligned}
q_{p}(t) &= \mu + \sum_{i=1}^2 \left( A_i \cos(i \omega_{day} t) + B_i \sin(i \omega_{day} t) \right) \\
&\quad \quad + \sum_{j=1}^2 \left( A'_j \cos(j \omega_{year} t) + B'_j \sin(j \omega_{year} t) \right) \\
&\quad \quad + \sum_{i=1}^2 \sum_{j=1}^2 \left( C_{ij} \cos(i \omega_{day} t)\cos(j \omega_{year} t) + D_{ij} \cos(i \omega_{day} t)\sin(j \omega_{year} t) \right) \\
&\quad \quad + \sum_{i=1}^2 \sum_{j=1}^2 \left( E_{ij} \sin(i \omega_{day} t)\cos(j \omega_{year} t) + F_{ij} \sin(i \omega_{day} t)\sin(j \omega_{year} t) \right)\\
\end{aligned}
\end{equation}

\end{document}